\documentclass[10pt,a4paper,twoside,english]{article}

\usepackage{cicis2}
\usepackage{blindtext}

\usepackage{multirow}

\begin{document}
\setcounter{page}{580}

\AddToShipoutPicture*{\BackgroundPic}

\title{Computation in Logic and Logic in
Computation\footnote{This paper is dedicated to Alan Turing, to commemorate
 the Turing Centenary Year 2012 -- his 100$^{\rm th}$ birthyear.}}
\author{
    Saeed Salehi\\
    \begin{small}
Department of Mathematical Sciences, University of Tabriz, 29 Bahman Blvd.,
51666--17766 Tabriz, Iran
    \end{small}\\
    \begin{small}
School of Mathematics, Institute for Research in Fundamental  Sciences (IPM),
19395--5746 Tehran, Iran
    \end{small}\\
    \begin{small}
http://saeedsalehi.ir/ \qquad root@saeedsalehi.ir \qquad
salehipour@tabrizu.ac.ir \end{small}
}
\date{}

\maketitle
\thispagestyle{empty}
\begin{cicisabstract}
The theory of addition in the domains of natural ($\mathbb{N}$), integer ($\mathbb{Z}$),
rational ($\mathbb{Q}$), real ($\mathbb{R}$) and complex ($\mathbb{C}$)
numbers is {\sf decidable}; so is the theory of multiplication in all those domains.
By G\"odel's Incompleteness Theorem the theory of addition and multiplication is
 {\sf undecidable} in the domains of $\mathbb{N}$, $\mathbb{Z}$ and $\mathbb{Q}$;
 though  Tarski proved that this theory is {\sf decidable} in the domains of
 $\mathbb{R}$ and $\mathbb{C}$. The theory of multiplication and order
  $\langle \cdot,\leqslant \rangle$ behaves differently in the above
  mentioned domains of numbers. By a theorem of Robinson, addition is definable by
  multiplication and order in the domain of natural numbers; thus the theory
    $\langle \mathbb{N}, \cdot,\leqslant\rangle$ is {\sf undecidable}. By a classical
     theorem in mathematical logic, addition is not definable in terms of multiplication
     and order in  $\mathbb{R}$. In this paper, we extend Robinson's theorem to the domain
      of integers ($\mathbb{Z}$)  by showing the definability of addition in
      $\langle \mathbb{Z},\cdot,\leqslant\rangle$; this implies that
        $\langle \mathbb{Z},\cdot,\leqslant\rangle$ is {\sf undecidable}.
        We also show the {\sf decidability} of  $\langle \mathbb{Q},\cdot,\leqslant\rangle$
         by the method of quantifier elimination. Whence, addition is not definable in
         $\langle \mathbb{Q},\cdot,\leqslant\rangle$.
\end{cicisabstract}

\begin{ciciskeywords}
Decidability; First-Order Logic; G\"odel's Incompleteness Theorems;  Church's Theorem;
 Presburger Arithmetic; Skolem Arithmetic; Quantifier Elimination.
\end{ciciskeywords}

\begin{cicismain}
\section{Introduction}
The question of the decidability of logical inference has triggered the beginning of
 computer science. Propositional Logic is decidable, since truth tables provide a finite
  semantics for it. Aristotle's Syllogism, or in modern terminology the first-order logic
   of unary predicates, is decidable, since it has the finite model property.  The notion
   of a Turing Machine was a successful outcome of the struggle to settle the question of
   the decidability of full First-Order Logic. It is now  known that the first-order logic
    is undecidable if it has a binary relation symbol or a binary function symbol
    (\cite{cdp}).  The additive theory of natural numbers $\langle\mathbb{N},+\rangle$
    was shown to be decidable by Presburger in 1929 (and by Skolem in 1930; see \cite{lnt}).
    The additive theories of integer, rational, real and complex numbers
    ($\langle\mathbb{Z},+\rangle$, $\langle\mathbb{Q},+\rangle$,
     $\langle\mathbb{R},+\rangle$ and $\langle\mathbb{C},+\rangle$) are
     decidable as well. The multiplicative theory of the natural numbers
     $\langle\mathbb{N},\cdot\rangle$ is also shown to be decidable by Skolem in 1930;
      the theories $\langle\mathbb{Z},\cdot\rangle$, $\langle\mathbb{Q},\cdot\rangle$,
       $\langle\mathbb{R},\cdot\rangle$ and $\langle\mathbb{C},\cdot\rangle$ are also
        decidable.

Then it was expected that the theory of addition and multiplication of natural numbers
would be decidable too; confirming Hilbert's Program. But the world was shocked in 1931
 by G\"odel's Incompleteness Theorem who showed that the theory
  $\langle\mathbb{N},+,\cdot\rangle$ is undecidable (see \cite{lnt}).
  The theory $\langle\mathbb{Z},+,\cdot\rangle$ is undecidable too,
  since $\mathbb{N}$ is definable in this structure: by Lagrange's
  Theorem $k\!\in\!\mathbb{N}\iff \exists \,a,b,c,d\!\in\!\mathbb{Z}\, (k=a^2+b^2+c^2+d^2)$.
   So is the theory $\langle\mathbb{Q},+,\cdot\rangle$ by Robinson's result \cite{robinson}
   which shows that $\mathbb{N}$ is definable in this structure too. However, Tarski showed
    that the theories $\langle\mathbb{R},+,\cdot\rangle$ and
    $\langle\mathbb{C},+,\cdot\rangle$ are decidable (\cite{marker}).
    It is worth mentioning that the order relation $\leqslant$ is definable by
    means of addition and multiplication in all the above domains of numbers.
    For example, the formulas $\exists z (z+x=y)$ and $\exists z (z^2+x=y)$ define
     the relation $x\leqslant y$ in the structures $\langle\mathbb{N},+,\cdot\rangle$
     and $\langle\mathbb{R},+,\cdot\rangle$ respectively. The theory of addition and
     order $\langle +,\leqslant\rangle$ is somehow weak, in all the above number domains,
     since it cannot define multiplication. The theory of multiplication and order
      $\langle\cdot,\leqslant\rangle$ has not been extensively studied; one reason is
       that addition is not definable in $\langle\mathbb{R},\cdot,\leqslant\rangle$,
       since the bijection $x\mapsto x^3$ of $\mathbb{R}$ preserves multiplication and
        order but does not preserve addition. Also it is known that addition is
        definable in $\langle\mathbb{N},\cdot,\leqslant\rangle$ by Tarski's identity
         (\cite{robinson}):
\centerline{$x+y=z \iff [x=y=z=0]\;\vee$}
\centerline{$[z\not=0\wedge \textsf{S}(z\cdot x)\cdot \textsf{S}(z\cdot y)=
\textsf{S}(z\cdot z\cdot \textsf{S}(x\cdot y))],$}
\\
where $\textsf{S}(u)$ is the successor of $u$, which is definable by
the order relation: $\textsf{S}(u)=v\!\iff\!\forall
w[u<w\leftrightarrow v\leqslant w]$. The symbol $u<v$ is a shorthand
for $u\leqslant v \wedge u\not=v$.

 The question of the decidability or undecidability of the structures
  $\langle\mathbb{Z},\cdot,\leqslant\rangle$ and $\langle\mathbb{Q},\cdot,\leqslant\rangle$
  are missing in the literature. In this paper, by modifying Tarski's identity we show
  that addition is definable in the structure
  $\langle\mathbb{Z},\cdot,\leqslant\rangle$;
  this implies the undecidability of $\langle\mathbb{Z},\cdot,\leqslant\rangle$.
  On the contrary, addition is not definable in $\langle\mathbb{Q},\cdot,\leqslant\rangle$;
  here we show a stronger result by the method of quantifier elimination:
  the theory $\langle\mathbb{Q},\cdot,\leqslant\rangle$ is decidable. Whence, by Robinson's
  above-mentioned result \cite{robinson}, addition cannot be defined in this structure.
  An interesting outlook of our results is that though $\langle +,\cdot\rangle$ puts
  the domains $\mathbb{N}$, $\mathbb{Z}$ and $\mathbb{Q}$ on the undecidable side,
  and the domains $\mathbb{R}$ and $\mathbb{C}$ on the decidable side,
  the language $\langle\cdot,\leqslant\rangle$ puts the domains $\mathbb{N}$
  and $\mathbb{Z}$ on the undecidable side, but $\mathbb{Q}$ and $\mathbb{R}$ on the
  decidable side.

\section{Multiplication and Order in $\mathbb{Z}$}
Tarski's identity $\textsf{S}(z\cdot x)\cdot \textsf{S}(z\cdot y)=\textsf{S}(z\cdot z\cdot
 \textsf{S}(x\cdot y))$ can define the formula $x+y=z$ in $\mathbb{Z}$ when $x+y\not=0$.
 The case $x+y=0$ was easily settled    in natural numbers: for any $x,y\in\mathbb{N}$ we
 have  $x+y=0\iff x=y=0$. But this does not hold in $\mathbb{Z}$, and so we have to treat
 this case differently. Our trick is to define the relation $x=-y$ in terms of multiplication
  and successor (which is definable by order): $x=-y\iff \textsf{S}(x)\cdot\textsf{S}(y)=
  \textsf{S}(x\cdot y)$. Thus, the following formula defines addition in terms of
  multiplication and order in $\mathbb{Z}$:
\centerline{$x+y=z \iff [z=0\wedge\textsf{S}(x)\cdot\textsf{S}(y)=
\textsf{S}(x\cdot y)]\;\vee$}
\centerline{$[z\not=0\wedge \textsf{S}(z\cdot x)\cdot \textsf{S}(z\cdot y)=
\textsf{S}(z\cdot z\cdot \textsf{S}(x\cdot y))].$}
\\
So, the theories $\langle\mathbb{Z},\cdot,\leqslant\rangle$ and
$\langle\mathbb{Z},+,\cdot\rangle$ are  interdefinable, and hence
 $\langle\mathbb{Z},\cdot,\leqslant\rangle$ is undecidable.

\section{Multiplication and Order in $\mathbb{Q}$}
Unlike the case of $\mathbb{Z}$, addition is not definable in the structure
$\langle\mathbb{Q},\cdot,\leqslant\rangle$. In fact, the theory of this structure
 is decidable. For showing that we use the method of quantifier elimination.
 First let us note that the language $\langle\cdot,\leqslant\rangle$ does not
 allow quantifier elimination for $\langle\mathbb{Q},\cdot,\leqslant\rangle$, since e.g.
 the formula $\exists y[x=y^2]$ is not equivalent to a quantifier-free formula. So,
 we restrict our attention to $\mathbb{Q}^+=\{r\in\mathbb{Q}\mid r>0\}$ and extend the
 language to $\mathcal{L}=\langle 0,1,\cdot, ^{-1},<,R_2,R_3,\ldots\rangle$, where $R_n$
 is interpreted as ``being the $n$th power of a rational"; or in other words
  $R_n(x)\equiv \exists y[x=y^n]$.

\noindent {\bf Theorem}.
The structure $\langle\mathbb{Q}^+,\mathcal{L}\rangle$ admits quantifier elimination.

We note that the above main theorem implies that the structure
$\langle\mathbb{Q},\mathcal{L}\rangle$ admits quantifier elimination
as well. It is enough to distinguish the signs: for any $x$, either
$-x>0$ or $x=0$ or $x>0$; so eliminating the quantifiers in each
case, will eliminate all of the quantifiers. Let us also note that
the quantifier-free formulas of $\mathcal{L}$ are decidable: for any
given rational number $r$ and any natural $n$ one can decide if $r$
is an $n$th power of (an-)other rational number or not. Thus,
quantifier elimination in $\langle\mathbb{Q},\mathcal{L}\rangle$
implies the decidability of the structure
$\langle\mathbb{Q},\mathcal{L}\rangle$, and hence $\langle\mathbb{Q},\cdot,\leqslant\rangle$.

\noindent The rest of the paper is devoted to proving the main
theorem. The folklore technique of quantifier elimination starts
from characterizing the terms and  atomic formulas, also eliminating
 negations, implications and universal quantifiers, and then
removing the disjunctions from the scopes of existential
quantifiers, which leaves the final case to be the existential
quantifier with the conjunction of some atomic (or negated atomic)
formulas. Removing this one existential quantifier implies the
ability to eliminate all the other quantifiers by induction. Let us
summarize the first steps:

For a variable $x$ and parameter ${a}$, all $\mathcal{L}-$terms are equal to $x^k{a}^l$ for
some $k,l\in\mathbb{Z}$. Atomic $\mathcal{L}-$formulas are in the form $u=v$ or $u<v$
or $R_n(u)$ for some terms $u,v$ and $n\geqslant 2$. Negated atomic $\mathcal{L}-$formulas
are thus $u\not=v$, $u\not<v$ and $\neg R_n(u)$; the formulas $u\not=v$ and $u\not<v$ are
 equivalent to $u<v \vee v<u$ and $u=v \vee v<u$ respectively. By de Morgan's laws we
 can assume that the negation appears only behind the atomic formulas of the form $R_n(u)$,
 and by the equivalences $A\rightarrow B \equiv \neg A\vee B$
 and $\forall x\varphi \equiv \neg\exists x\neg \varphi$, we can assume that
 the implication symbol and universal quantifier do not appear in the formula
 (whose quantifiers are to eliminated). Finally,
  the equivalence $\exists x (\varphi\vee\psi) \equiv \exists x\varphi \vee \exists x\psi$
  leaves us with the elementary formulas of the form $\exists x (\bigwedge_i \theta_i)$
  where each $\theta_i$ is in the form $(x^\alpha=v)$ or $(r<x^\beta)$ or $(x^\gamma<s)$
  or $R_n(tx^\delta)$ or $\neg R_m(ux^\epsilon)$ for some
   $\alpha,\beta,\gamma,\delta,\epsilon\in\mathbb{N}$ and $\mathcal{L}-$terms $r,s,t,u,v$.
    Whence, it suffices to show that the $\mathcal{L}-$formula
  $\exists x \big[\bigwedge_h(x^{\alpha_h}=v_h)\wedge \bigwedge_i(r_i<x^{\beta_i})
  \wedge\bigwedge_j(x^{\gamma_j}<s_j)\wedge \bigwedge_k(R_{n_k}(t_k\cdot x^{\delta_k}))
  \wedge\bigwedge_l(\neg R_{m_l}(u_l\cdot x^{\epsilon_l}))\big]$ is equivalent to another
  $\mathcal{L}-$formula in which $x$ (and so $\exists x$) does not appear.
  This will finish the proof.

Here comes the next steps of quantifier elimination.  The powers of
$x$ can be unified: let $p$  be the least common multiplier of the
$\alpha_h$'s, $\beta_i$'s, $\gamma_j$'s, $\delta_k$'s and
$\epsilon_l$'s.  From the
$\langle\mathbb{Q}^+,\mathcal{L}\rangle-$equivalences $a\!=\!b
\leftrightarrow a^q\!=\!b^q$, $a\!<\!b\leftrightarrow a^q\!<\!b^q$
and $R_n(a)\leftrightarrow R_{nq}(a^q)$,  we infer that the above
formula can be re-written equivalently as
\\
$\exists x \big[\bigwedge_h(x^{p}\!=\!v_h)\wedge \bigwedge_i(r_i\!<\!x^{p})
\wedge\bigwedge_j(x^{p}\!<\!s_j)\wedge$
\\
$\bigwedge_k(R_{n_k}(t_k\!\cdot\!x^{p}))\wedge\bigwedge_l(\neg R_{m_l}(u_l\!\cdot\!x^{p}))
\big]$
\\
for possibly   new  $v_h$'s, $r_i$'s, $s_j$'s, $n_k$'s, $t_k$'s,
$m_l$'s and $u_l$'s. This formula is in turn  equivalent to
 \\
 $\exists y\big[\bigwedge_h(y\!=\!v_h)\wedge \bigwedge_i(r_i\!<\!y)
 \wedge\bigwedge_j(y\!<\!s_j)\wedge$
 \\
$\bigwedge_k(R_{n_k}(t_k\!\cdot\!y))\wedge\bigwedge_l(\neg R_{m_l}(u_l\!\cdot\!y))
\wedge R_p(y)\big]$
  \\
  (with the substitution $y=x^p$).
  Thus it suffices to show that the following formula
   \\
   $\exists x \big[\bigwedge_h(x=v_h)\wedge \bigwedge_i(r_i<x)\wedge\bigwedge_j(x<s_j)\wedge$
   \\
    $\bigwedge_k(R_{n_k}(t_k\cdot x))\wedge\bigwedge_l(\neg R_{m_l}(u_l\cdot x))\big]$
   \\
   is equivalent to a quantifier-free formula.
 If the conjunction $\bigwedge_h(x=v_h)$ is not empty, then the above formula is
 equivalent to the quantifier-free formula
 \\
 $[\bigwedge_h(v_0=v_h)\wedge \bigwedge_i(r_i<v_0)\wedge\bigwedge_j(v_0<s_j)\wedge$
 \\
  $\bigwedge_k(R_{n_k}(t_k\cdot v_0))\wedge\bigwedge_l(\neg R_{m_l}(u_l\cdot v_0))\big]$
  \\
  for some term $v_0$. So, let us assume that the conjunction $\bigwedge_h(x=v_h)$ is
  empty, and thus we are to eliminate the quantifier of the formula
  \\
$\exists x \big[\bigwedge_i(r_i<x)\wedge\bigwedge_j(x<s_j)\wedge
\bigwedge_k(R_{n_k}(t_k\cdot x))\wedge\bigwedge_l(\neg R_{m_l}(u_l\cdot x))\big]$.

The formula $\exists x \big[\bigwedge_i(r_i\!<\!x)\wedge\bigwedge_j(x\!<\!s_j)\big]$ is
$\langle\mathbb{Q}^+,\mathcal{L}\rangle-$equivalent to (the quantifier-free formula)
$\bigwedge_{i,j}(r_i\!<\!s_j)$ (that is $\max_i\{r_i\}\!<\!\min_j\{s_j\}$),
since $\langle\mathbb{Q}^+,<\rangle$ is dense.

For the formula $\exists x \bigwedge_k R_{n_k}(t_k\cdot x)$, let $\textsf{p}$ be a
 prime number, and put $t_k'$ be the greatest number such that $\textsf{p}^{t_k'}$
 divides $t_k$; similarly $x'$ is the greatest number such that $\textsf{p}^{x'}$
 divides $x$. Then $\bigwedge_k R_{n_k}(t_k\cdot x)$ is equivalent to
 $\forall\textsf{p}\bigwedge_k[t_k'+x'\equiv_{n_k}0]$. By a generalized form of
the Chinese Remainder Theorem (\cite{lnt}) the existence of such an
$x'$ is equivalent
  to $\bigwedge_{\kappa\not=\lambda}\;
  t_\kappa'\equiv_{(n_\kappa,n_\lambda)}t_\lambda'$;
  here $(a,b)$ is the greatest common divisor of $a$ and $b$.
That is equivalent to
  $\bigwedge_{\kappa\not=\lambda}R_{(n_\kappa,n_\lambda)}(t_\kappa\cdot t_\lambda^{-1})$.
  We further note that in case of
  $\bigwedge_{\kappa\not=\lambda}\;t_\kappa'\equiv_{(n_\kappa,n_\lambda)}t_\lambda'$
   there are infinitely many solutions for  $\bigwedge_k[t_k'+x'\equiv_{n_k}0]$ which
    are in the form  $x'=Ny' - \sum_k \nu_kt_k'$ for some fixed integers $N$ and
     $\nu_k$'s; $y'$ is arbitrary. In fact $N$ is the least common multiplier of $n_k$'s, and
$\nu_k$'s are $\nu_k=c_k\!\cdot\!N/n_k$ where $\sum_k c_kN/n_k=1$; the existence of 
$c_k$'s follows from the fact that the greatest common divisor of
$(N/n_k)$'s is $1$. Moreover, the solution $x'$ is unique up to the
module $N$.
     So, if there exists some $x\in\mathbb{Q}^+$ which
     satisfies $\bigwedge_k R_{n_k}(t_k\cdot x)$ for some $t_k\in\mathbb{Q}^+$,
     then it must be of the form $x=\gamma^N\cdot\prod_k(t_k)^{-\nu_k}$ for some (arbitrary)
     $\gamma\in\mathbb{Q}^+$.

Thus, the formula $\exists x \big[\bigwedge_i(r_i<x)\wedge\bigwedge_j(x<s_j)\wedge
\bigwedge_k(R_{n_k}(t_k\cdot x))\big]$ is equivalent to (the quantifier-free formula)
 $\bigwedge_{i,j}(r_i\!<\!s_j)\wedge\bigwedge_{\kappa\not=\lambda}R_{(n_\kappa,n_\lambda)}
 (t_\kappa\cdot t_\lambda^{-1})$, since the solution $x=\gamma^N\cdot\prod_k(t_k)^{-\nu_k}$
 for $\bigwedge_k(R_{n_k}(t_k\cdot x))$ can be chosen to satisfy
 $\max_i\{r_i\}<x<\min_j\{s_j\}$: choose a rational number $\gamma\in\mathbb{Q}^+$
 between the positive real numbers
 $\alpha=\big(\max_i\{r_i\}\cdot(\prod_k(t_k)^{\nu_k})\big)^{1/N}$ and
 $\beta=\big(\min_j\{s_j\}\cdot(\prod_k(t_k)^{\nu_k})\big)^{1/N}$.
 Since the set $\mathbb{Q}$ is dense in $\mathbb{R}$, there exists such a rational
 number $\gamma$.
 Then  $x=\gamma^N\cdot\prod_k(t_k)^{-\nu_k}$ is the desired solution.

Finally, we show that the formula
\\
$\exists x \big[\bigwedge_i(r_i<x)\wedge\bigwedge_j(x<s_j)\wedge
\bigwedge_k(R_{n_k}(t_k\cdot x))\wedge$
\\
$\bigwedge_l(\neg R_{m_l}(u_l\cdot x))\big]$
\\
is equivalent to  the following quantifier-free formula
\\
$\bigwedge_{i,j}(r_i<s_j)\wedge\bigwedge_{\kappa\not=\lambda}R_{(n_\kappa,n_\lambda)}
(t_\kappa\!\cdot\!t_\lambda^{-1})\wedge$
\\
$\bigwedge_{l:m_l\mid N}(\neg R_{m_l}(u_l\!\cdot\!t))$,
\\
where $N$ is the least common multiplier of $n_k$'s, and
$t=\prod_\kappa (t_\kappa)^{-\nu_\kappa}$ in which $\nu_k=c_kN/n_k$'s satisfy
$\sum_k c_kN/n_k=1$.

If for some $x\in\mathbb{Q}^+$,
$\bigwedge_i(r_i<x)\wedge\bigwedge_j(x<s_j)\wedge
\bigwedge_k(R_{n_k}(t_k\cdot x))\wedge\bigwedge_l(\neg
R_{m_l}(u_l\cdot x))$ holds, then clearly $\bigwedge_{i,j}(r_i<s_j)$
is true, and it can be easily seen that we also have
$\bigwedge_{\kappa\not=\lambda}R_{(n_\kappa,n_\lambda)}
(t_\kappa\!\cdot\!t_\lambda^{-1})$. Assume $m_l\mid N$; we show that
$\neg R_{m_l}(u_l\cdot t)$. Note that there exists some $\gamma$
such that $x=\gamma^N\cdot t$. Now if $R_{m_l}(u_l\cdot t)$, then
$u_l\cdot x=\gamma^N\cdot u_l\cdot t$, and so by $m_l\mid N$ we have
$R_{m_l}(u_l\cdot x)$ which contradicts the assumption
$\bigwedge_l(\neg R_{m_l}(u_l\cdot x))$. Whence,
$\bigwedge_{l:m_l\mid
 N}(\neg R_{m_l}(u_l\!\cdot\!t))$ holds.

 Conversely, if we have
$\bigwedge_{i,j}(r_i<s_j)\wedge\bigwedge_{\kappa\not=\lambda}R_{(n_\kappa,n_\lambda)}
(t_\kappa\cdot t_\lambda^{-1})\wedge\bigwedge_{l:m_l\mid N}(\neg
R_{m_l}(u_l\cdot t))$, then by the above arguments there exist some
positive real numbers $\alpha < \beta$ such that for any rational
$\gamma$ with $\alpha<\gamma<\beta$,  the number $z=\gamma^N\cdot t$
satisfies the formula
$\bigwedge_i(r_i<z)\wedge\bigwedge_j(z<s_j)\wedge
\bigwedge_k(R_{n_k}(t_k\cdot z))$ where $N$ and $t$ are as above.
Let $\textsf{P}$ be a sufficiently large prime number which does not
divide any of the numerators or denominators  of (the reduced
fractions of) $t_k$'s or $u_l$'s. Let $M=\prod_lm_l$ and let
$\delta$  be a positive rational number such that
$(\alpha/\textsf{P})^{1/M}<\delta<(\beta/\textsf{P})^{1/M}$. We
show that $x=\textsf{P}^N\cdot \delta^{N\cdot M}\cdot t$ satisfies
$\bigwedge_l\neg R_{m_l}(u_l\cdot x)$. Note that since
$\alpha<\textsf{P}\cdot\delta^M<\beta$ we already have
$\bigwedge_i(r_i<x)\wedge\bigwedge_j(x<s_j)\wedge
\bigwedge_k(R_{n_k}(t_k\cdot x))$. For showing $\neg
R_{m_l}(u_l\cdot x)$ we distinguish two cases. (1) If $m_l\mid N$
then $R_{m_l}(u_l\cdot x)\equiv R_{m_l}(u_l\cdot \textsf{P}^N\cdot
\delta^{N\cdot M}\cdot t)$ implies $R_{m_l}(u_l\cdot t)$
contradicting the assumption $\bigwedge_{l:m_l\mid N}(\neg
R_{m_l}(u_l\cdot t))$; thus $\neg R_{m_l}(u_l\cdot x)$. (2) If
$\neg(m_l\mid N)$, then $R_{m_l}(u_l\cdot x)$ or equivalently
$R_{m_l}(u_l\cdot \textsf{P}^N\cdot \delta^{N\cdot M}\cdot t)$
implies $R_{m_l}(u_l\cdot t\cdot \textsf{P}^N)$ since $m_l\mid M$.
Since $\textsf{P}$  does not divide any of the numerators or
denominators of (the reduced fractions of) $u_l$'s or $t$ ($t_k$'s), then we must have $R_{m_l}(\textsf{P}^N)$ which holds if
and only if $m_l\mid N$; this contradicts our assumption
$\neg(m_l\mid N)$. Thus $\neg R_{m_l}(u_l\cdot x)$. Whence, all in
all we showed that $\bigwedge_l\neg R_{m_l}(u_l\cdot x)$ holds.
 \hfill  Q.E.D

\paragraph{Acknowledgements} This research was partially supported
by grant  No.~90030053 of the Institute for Research in Fundamental Sciences (IPM), Tehran.

\end{cicismain}

\begin{center}
\begin{tabular}{|c|}
\hline
\\
{\sc Saeed Salehi}, ``Computation in Logic and Logic in Computation", Invited Paper in: \\
{Bahram B. Sadeghi} (editor), {\it Proceedings of the Third International Conference on } \\
{\it Contemporary Issues in Computer and Information Sciences} (CICIS 2012), \\
Institute for Advanced Studies in Basic Sciences, Gavazangh, Zanjan, Iran, \\
Brown Walker Press 2012, USA, ISBN 9781612336237 (624 pages), pp. 580--583.\\
\url{http://www.universal-publishers.com/book.php?method=ISBN&book=161233623X}\\
\\
\hline
\end{tabular}
\end{center}


\section*{References}
\begin{biblist}

\bib{cdp}{book}{
title={The Classical Decision Problem},
author={B\"orger,E.},
author={Gr\"adel,E.},
author={Gurevich,Y.},
date={2001},
publisher={Springer-Verlag, Berlin},
address={}
}

\bib{robinson}{article}{
title={Definability and Decision Problems in Arithmetic},
subtitle={},
author={Robinson,J.},
author={ },
author={ },
journal={The Journal of Symbolic Logic},
volume={14},
date={1949},
pages={98--114}
}

\bib{marker}{book}{
title={Model Theory: An Introduction},
author={Marker,D.},
date={2002},
publisher={Springer-Verlag, Berlin},
address={}
}



\bib{lnt}{book}{
title={Logical Number Theory I: An Introduction},
author={Smory\'nski,C.},
date={1991},
publisher={Springer-Verlag, Berlin},
address={}
}



\end{biblist}
\end{document}